\input amstex.tex
\input amsppt.sty
\input xy.tex
\xyoption{all}
\magnification1200
\NoPageNumbers
\nologo
\TagsOnRight
\voffset-30pt

\def\bs{{\bar s}}
\def\C{{\Cal C}}
\def\CC{{\Bbb C}}
\def\crc{{\raise.24ex\hbox{$\sssize\kern.1em\circ\kern.1em$}}}
\def\de{\delta}

\def\deg{{\operatorname{deg}\ts\ts}}
\def\di{{\operatorname{dim}}}
\def\End{{\operatorname{End}}}
\def\fm{F}
\def\G{{\Cal G}}
\def\ga{\gamma}
\def\Ga{\Gamma}
\def\gr{{\operatorname{gr}}\hskip1pt}
\def\Hen{{\text{\it H}\hskip1pt(n)}}
\def\Hom{{\operatorname{Hom}}}

\def\L{{\Cal L}}
\def\la{\lambda}
\def\La{\Lambda}
\def\Lap{{\La^\prime}}
\def\Lapp{{\La^{\prime\prime}}}
\def\np{{\par\noindent}}
\def\O{{\Cal O}}
\def\Om{\Omega}
\def\Omp{{\Om^\prime}}
\def\Ompp{{\Om^{\ts\prime\prime}}}
\def\ot{\otimes}
\def\ol{\overline}
\def\ph{\varphi}
\def\pr{{^{\ts\prime}}}
\def\prpr{{^{\prime\prime}}}
\def\RR{{\Bbb R}}
\def\si{{\operatorname{sign}}}
\def\sq{{\ $\square$}}
\def\ts{\thinspace}
\def\tts{\hskip.5pt}

\def\ty{{u}}
\def\Wen{{\text{\it A}\hskip.5pt(n,N)}}

\def\gad#1{\global\advance#1 by1}
\newcount\Sno
\def\reset{\gad\Sno\Fno=0}
\newcount\Fno
\def\Tag#1{\expandafter\ifx\csname#1\endcsname\relax\gad\Fno
\expandafter\xdef\csname#1\endcsname{\the\Sno.\the\Fno}\fi\tag{\[#1]}}
\def\[#1]{\csname#1\endcsname}

\def\Lac{{\mathchoice
 {\textstyle\La^{\kern-2pt\circ}}{\textstyle\La^{\kern-2pt\circ}}
 {\scriptstyle\La^{\kern-2pt\raise-.5pt\hbox{$\scriptscriptstyle\circ$}}}{ch}}}

\font\bigbf=cmbx10 scaled 1200
\font\smallrm=cmr8 scaled 1000

\centerline{\bigbf Young's Orthogonal Form for}
\smallskip
\centerline{\bigbf Brauer's Centralizer Algebra}
\vskip22pt
\centerline{\smc Maxim Nazarov}
\vskip10pt
\topmatter
\abstract
\nofrills
We consider the semi-simple algebra which arises as the centralizer
of a tensor power of the fundamental representation of the orthogonal group.
There is a canonical basis in every irreducible representation of this
algebra; it is an analogue of the Young basis in an irreducible representation
of the symmetric group. We evaluate the action of the generators of this
centalizer algebra in the canonical basis. We use this result to define
an analogue of the degenerate affine Hecke algebra of the general linear group.
\endabstract
\endtopmatter

\heading
\vskip-10pt
Introduction
\endheading

\np
Let $G$ be one of the classical groups
$GL(N,\CC)$, $O(N,\CC)$, $Sp(N,\CC)$
acting on the vector space $U=\CC^N$.
The question of how the $n$-th tensor power of the representation
$U$ decomposes into irreducible summands
leads to studying the centralizer $C(n,N)$ in 
$\End(U)^{\ot n}$ of the image of the group $G$. 
By the definition of the algebra $C(n,N)$ we have the
ascending chain of subalgebras
$$
C(1,N)\subset C(2,N)\subset\ldots\subset C(n,N).
$$ 
Moreover, for the classical group $G$
any irreducible representation of
$C(n,N)$ appears at most once in the restriction of an irreducible
representation of $C(n+1,N)$. Therefore a canonical basis
exists in any irreducible representation $V$ of $C(n,N)$. Its
vectors are the eigenvectors for the subalgebra $X(n,N)$ in $C(n,N)$
generated by all the central elements in the members of the
above chain.

For the group $G=GL(N,\CC)$ the centralizer $C(n,N)$
is generated by the permutational action of the symmetric group
$S(n)$ in $U^{\ot n}$. The action of $S(n)$ on the vectors
of the canonical basis in $V$
was described for the first time by A.\ts Young\ts\nolinebreak[Y].
G.\ts Murphy [Mp] rederived the formulas from [Y] by using the
properties of the subalgebra $X(n,N)$.

Now let $G$ be the orthogonal group $O(N,\CC)$.
To describe the corresponding centralizer algebra $C(n,N)$
explicitly,
R.\ts Brauer [Br] introduced a certain complex associative algebra
$B(n,N)$ along with a homomorphism onto $C(n,N)$.
This homomorphism is injective if and only if $N\geqslant n$.
There is also a chain of subalgebras
$$
B(1,N)\subset B(2,N)\subset\ldots\subset B(n,N).
$$ 
The group algebra $\CC[S(n)]$ is contained in $B(n,N)$
as a subalgebra.
The structure of the algebra $B(n,N)$ was investigated by
P.\ts Hanlon and D.\ts Wales; see [HW] and references therein. 
In the present article we will also work with
$B(n,N)$ and regard $V$ as a representation of the
latter algebra. 

For $N\geqslant n$ an explicit description of the action of the
algebra $B(n,N)$ on the vectors of the canonical basis in $V$
was given by J.\ts Murakami [Mk].
His description was based on the results of [JMO]. 
In the present article for any $N$ we give a new description of
this action based entirely on the properties of the
subalgebra $X(n,N)$ in $C(n,N)$. The case $G=Sp(N,\CC)$
is very similar and will be considered elsewhere.  


In Section 2 we introduce a remarkable family of pairwise
commuting elements $x_1,\dots,x_n$ of the algebra $B(n,N)$.
For every $n$ the element $x_{n+1}$ belongs to the
centralizer of the subalgebra $B(n,N)$ in $B(n+1,N)$.
The elements $x_1,\dots,x_n$ are the analogues of the
pairwise commuting elements of $\CC[S(n)]$
which were used in [Ju,Mu].
Their images in $C(n,N)$ belong to the subalgebra $X(n,N)$.
The vectors of the canonical basis in $V$ are eigenvectors 
of the elements $x_1,\dots,x_n$ and we evaluate the respective
eigenvalues; see Theorem 2.6.

There is a natural projection map $B(n+1,N)\to B(n,N)$
commuting with both left and right multiplication by the
elements from $B(n,N)$; this map has been already used by H.\ts Wenzl 
in [W].
The images of powers of the element
$x_{n+1}$ 
with
respect
to
this
map are
certain central elements of the algebra $B(n,N)$. We evaluate
the eigenvalues of these central elements in every
irreducible representation $V$; see Theorem 3.9.

The algebra $B(n,N)$ comes with a family
of generators $s_1,\dots,s_{n-1}$;
$\bs_1,\dots,\bs_{n-1}$.
The elements $s_1,\dots,s_{n-1}$ are the standard
generators of the symmetric group $S(n)$.
Moreover, the quotient of the algebra
$B(n,N)$ with respect to the ideal
generated by $\bs_1,\dots,\bs_{n-1}$ is isomorphic to $\CC[S(n)]$.
We point out certain relations between the elements
$x_1,\dots,x_n$ and the  generators of $B(n,N)$;
see Proposition 2.3. By using Proposition 2.3 and
Theorems 2.6,\ts 3.9
we describe the action of these generators on
the vectors
of the canonical basis in every representation $V$.
For the representations which factorize through $\CC[S(n)]$
our formulas coincide with those from [Y].
 
In Section 4 we 
use the results of Sections 2 and 3 as a motivation to 
introduce a new algebra. This is an analogue of
the degenerate affine Hecke algebra $\Hen$ from [C1,C2] and [D].
We denote the new algebra by $\Wen$ and call it the
affine Brauer algebra. The algebra $\Hen$ is a quotient
of $\Wen$; see Corollary~4.9.
For each $m=0,1,2,\ldots$ the
algebra $\Wen$ admits a homomorphism to the centralizer 
of the subalgebra $B(m,N)$ in $B(m+n,N)$. The kernels
of all these homomorphisms have the zero intersection;
see Theorem 4.7. We use these homomorphisms to
construct a linear basis in the algebra $\Wen$; see Theorem~4.6.

I am very grateful to
D.\ts E.\ts Evans,
A.\ts O.\ts Morris,
G.\ts I.\ts Olshanski
and
A.\ts M.\ts Vershik
for numerous discussions. I am also very grateful to all my colleagues
at the Institut Gaspard Monge, Universit\'e Marne-la-Vall\'ee for
their generous hospitality.

\headline={\smallrm
\ifodd\pageno
Brauer's Centralizer Algebras\hfill\number\pageno
\else
\number\pageno\hfill Maxim Nazarov
\fi}

\heading
1. Brauer Centralizer Algebra
\endheading\reset

\np
Let $n$ be a positive integer and $N$ be an arbitrary complex parameter.
Denote by $\G(n)$ be the set of all graphs with $2n$ vertices and
$n$ edges such that each vertex is incident with an edge. We will
enumerate the vertices by $1,\ldots,n,\bar1,\ldots,\bar n$.
In other words, $\G(n)$ consists of all partitions
of the set $\{1,\ldots,n,\bar1,\ldots,\bar n\}$ into pairs. 
We will define the {\it Brauer algebra} $B(n,N)$ as an associative
algebra over $\CC$ with the basic elements $b(\ga),\ts\ga\in\G(n)$.

To describe the product $b(\ga)\ts b(\ga\pr)$ in $B(n,N)$ consider
the graph obtained by identifying the vertices $\bar1,\dots,\bar n$ of
$\ga$ with the vertices $1,\dots,n$ of $\ga\pr$ respectively.
Let $q$ be the number of loops in this graph. Remove
all the loops and replace the remaining connected components
by single edges, retaining the numbers of the terminal vertices.
Denote by $\ga\crc\ga\pr$ the resulting graph, then by definition
$$
b(\ga)\ts b(\ga\pr)=N^{\ts q}\cdot b(\ga\crc\ga\pr).
\Tag{1.0}
$$      

\newpage

Evidently, the dimension of $B(n,N)$ is equal to
$1\cdot3\cdot5\cdot\ldots\cdot(2n-1)$. The algebra $B(n,N)$ contains the
group algebra of the symmetric group $S(n)$; one can identify an element
$s$ of $S(n)$ with $b(\ga)$ where the edges of $\ga$ are
$\{s(1),\bar1\},\dots,\{s(n),\bar n\}$.

An edge of the form $\{k,\bar k\}$ will be called {\it vertical}.
We will regard $B(n-1,N)$ as a subalgebra of $B(n,N)$ with the basic
elements $b(\ga)$ where $\ga$ contains the vertical edge $\{n,\bar n\}$.   
Along with a transposition $(k,l)$ in $S(n)$ we will consider the element
$\ol{(k,l)}=b(\ga)$ of $B(n,N)$ where the only non-vertical edges
of $\ga$ are $\{k,l\}$ and $\{\bar k,\bar l\}$. 

We will sometimes write $s_k$ and $\bs_k$ instead of $(k,k+1)$ and
$\ol{(k,k+1)}$ respectively. The elements
$s_1,\dots,s_{n-1};\bs_1,\dots,\bs_{n-1}$ generate the algebra  
$B(n,N)$. One can directly verify the following relations for these
elements:
$$
\gather
s_k^2=1;
\quad
\bs_k^{\ts2}=N\ts\bs_k;
\quad
s_k\ts\bs_k=\bs_k\ts s_k=\bs_k;
\Tag{1.1}
\\
s_k\ts s_{k+1}\ts s_k=s_{k+1}\ts s_k\ts s_{k+1};
\quad
\bs_k\ts\bs_{k+1}\ts\bs_k=\bs_k;
\quad
\bs_{k+1}\ts\bs_k\ts\bs_{k+1}=\bs_{k+1};
\Tag{1.2}
\\
s_k\ts\bs_{k+1}\ts\bs_k=s_{k+1}\ts\bs_k;
\quad
\bs_{k+1}\ts\bs_k\ts s_{k+1}=\bs_{k+1}\ts s_k;
\Tag{1.3}
\\
s_k\ts s_{l}=s_{l}\ts s_k,
\quad
\bs_k\ts s_{l}=s_{l}\ts\bs_k,
\quad
\bs_k\ts\bs_{l}=\bs_{l}\ts\bs_k,
\quad\ |k-l|>1.
\Tag{1.4}
\endgather
$$

\proclaim{Proposition 1.1}
The relations {\rm (\[1.1])} to {\rm (\[1.4])} are defining
relations for $B(n,N)$.
\endproclaim

\np
For the proof of this proposition see [BW, Section 5\ts ].
Now suppose that $N$ is a positive integer. Consider the $n$-th
tensor power of the
representation $U=\CC^N$
of the orthogonal group $G=O(N,\CC)$. Let $u(1),\dots,u(N)$ be the
standard orthogonal basis in $U$; denote by
$u(i_1\ldots i_n)$ the vector $u(i_1)\ot\ldots\ot u(i_n)$ in $U^{\ot n}$.
Consider the centralizer  algebra
$C(n,N)=\End_{\ts G}\bigl(U^{\ot n}\bigr)$.

\proclaim{Proposition 1.2}
a) There is a homomorphism
$B(n,N)\rightarrow C(n,N)$
where the actions of $(k,l)$ and $\ol{(k,l)}$ in $U^{\ot n}$ for $k<l$
are defined by
$$
\align
(k,l)\cdot u(i_1\ldots i_k\ldots i_l\ldots i_n)
&=u(i_1\ldots i_l\ldots i_k\ldots i_n),
\Tag{1.5}
\\
\ol{(k,l)}\cdot u(i_1\ldots i_k\ldots i_l\ldots i_n)
&=\de(i_k\ts i_l)\cdot\sum_{i=1}^N\ts
u(i_1\ldots i\ldots i\ldots i_n).
\endalign
$$

b) This homomorphism is surjective for any positive integer $N$.

c) This homomorphism is injective if and only if $N\geqslant n$.
\endproclaim

\demo{Proof}
The actions of the elements $(k,l)$ and
$\ol{(k,l)}$ in $U^{\ot n}$ evidently commute with the action of
the orthogonal group $G$. The parts a) and b) are results
of [Br,\ts Section~5]. The part c) follows from
[B2,\ts Theorem 7A]\sq
\enddemo

\np
The algebra $C(n,N)$ is semisimple by its definition;
the irreducible representations of $C(n,N)$
are parametrized [Wy,\ts Theorem 5.7.F] by Young diagrams with
at most $N$ boxes in the first two columns and with $n-2r$ boxes
altogether where $r=0,1,\dots,[\ts{n}/2\ts]$. Denote the set of all
such diagrams by $\Cal O(n,N)$.
Let $V(\la,n)$ be the representation of $C(n,N)$
corresponding to a diagram $\la\in\O(n,N)$. The next proposition
is contained in [L,\ts Theorem I]; see also [Ki, Section 3].  

\proclaim{Proposition 1.3}
The restriction of $V(\la,n)$ to $C(n-1,N)$ decomposes
into the direct sum $\underset\mu\to\oplus\ts V(\mu,n-1)$ where $\mu$
ranges over all the diagrams $\mu\in\O(n-1,N)$ obtained from $\la$ by
removing or adding a box.
\endproclaim

\proclaim{Corollary 1.4}
Each irreducible representation of $C(n-1,N)$ appears at most once in the
restriction onto $C(n-1,N)$ of an irreducible representation of $C(n,N)$.
\endproclaim 

\heading
2. Jucys-Murphy Elements for $B(n,N)$
\endheading\reset

\np
By definition for any complex parameter $N$ we have the chain of
subalgebras
$$
B(1,N)\subset B(2,N)\subset\ldots\subset B(n,N).
\Tag{2.0}
$$
In this section we will introduce a remarkable family of pairwise
commuting elements in $B(n,N)$ corresponding to this chain; cf.
[Ju,\ts Mu].
For every $k=1,\dots,n$ consider the element of $B(k,N)$
$$
x_k=\frac{N-1}2+\sum_{l=1}^{k-1}\ts
(k,l)-\ol{(k,l)}.
\Tag{2.1}
$$

\proclaim{Lemma 2.1}
The element $x_k$ commutes with all the elements of $B(k-1,N)$.
\endproclaim

\demo{Proof}
The right hand side of (\[2.1]) is symmetric in $l=1,\dots,k-1$. 
Therefore $x_k$ commutes with any element $s$ of $S(n-1)$. To complete
the proof it suffices to check that $x_k$ commutes with $\bs_{k-2}$.
The commutator
$
\bigl[\ts\bs_{k-2}\ts,\ts x_k\ts\bigr]
$
equals
$$
\bigl[\ts\ol{(k-2,k-1)}\ts,\ts
(k-2,k)-\ol{(k-2,k)}+(k-1,k)-\ol{(k-1,k)}\ts\ts\bigr].
$$
The latter commutator vanishes because 
$$
\align
\ol{(k-2,k-1)}\cdot(k-2,k)=&\ts\ol{(k-2,k-1)}\cdot\ts\ol{(k-1,k)}\ts,
\\
\ol{(k-2,k-1)}\cdot\ol{(k-2,k)}=&\ts\ol{(k-2,k-1)}\cdot(k-1,k)\ts,
\\
(k-2,k)\cdot\ol{(k-2,k-1)}=&\ts\ol{(k-1,k)}\cdot\ol{(k-2,k-1)}\ts
\\
\ol{(k-2,k)}\cdot\ol{(k-2,k-1)}=&\ts(k-1,k)\cdot\ts\ol{(k-2,k-1)}\ts.
\endalign
$$
\nopagebreak
The last four equalities are verified directly by the definition 
(\[1.0])\sq
\enddemo

\proclaim{Corollary 2.2}
The elements $x_1,\ldots,x_n$ of $B(n,N)$ pairwise commute.
\endproclaim

\proclaim{Proposition 2.3}
The following relations hold in the algebra $B(n,N)$:
$$
\align
s_k\ts x_l=x_l\ts s_k,
&\quad
\bs_k\ts x_l=x_l\ts\bs_k;
\qquad
l\neq k,k+1;
\Tag{2.2}
\\
s_k\ts x_k-x_{k+1}\ts s_k=\bs_k-1,
&\quad
s_k\ts x_{k+1}-x_k\ts s_k=1-\bs_k;
\Tag{2.3}
\\
\bs_k\ts(x_k+x_{k+1})=0,
&\quad
(x_k+x_{k+1})\ts\bs_k=0.
\Tag{2.4}
\endalign
$$
\endproclaim

\demo{Proof}
The relations (\[2.2]) for $l>k+1$ follow from Lemma 2.1 while
those for $l<k$ follow directly from the definition (\[2.1]). Also
by this definition we have the equality
$$
x_{k+1}-s_k\ts x_k\ts s_k=s_k-\bs_k
$$
which implies the relations (\[2.3]). Again using (\[2.1]) we
obtain for any $l=1,\dots,k-1$ the equalities
$$
\gather
\ol{(k,k+1)}\cdot(k,l)=\ol{(k,k+1)}\cdot\ol{(k+1,l)},
\\
\ol{(k,k+1)}\cdot\ol{(k,l)}=\ol{(k,k+1)}\cdot(k+1,l)
\endgather
$$
which together with 
$$
\ol{(k,k+1)}\cdot(k,k+1)=\ol{(k,k+1)},
\quad
\ol{(k,k+1)}^{\ts2}=N\cdot\ol{(k,k+1)}
$$
imply the first relation in (\[2.4]). The proof of the second relation is
quite similar\sq
\enddemo

\proclaim{Corollary 2.4}
The elements $x_1^i+\ldots+x_n^i$ with $i=1,3,\ldots$
are central in $B(n,N)$.
\endproclaim 

\demo{Proof}
For any $i=1,2,3,\ldots$ the relations (\[2.3]) imply that 
$$
\gather
s_k\ts x_k^i=x_{k+1}^i\ts s_k
+\sum_{j=1}^i\ts\ts x_{k+1}^{j-1}\ts(\bs_k-1)\ts x_k^{i-j},
\Tag{2.425}
\\
s_k\ts x_{k+1}^i=x_{k}^i\ts s_k
-\sum_{j=1}^i\ts\ts x_{k}^{j-1}\ts(\bs_k-1)\ts x_{k+1}^{i-j}.
\endgather
$$
Combining (\[2.4]) with (\[2.425]) gives for odd $i$ the equalities
$$
\align
\bigl[s_k\ts,\ts x_k^i+x_{k+1}^i\ts\bigr]
&=
\sum_{j=1}^i\ts\ts
\bigl(\ts x_{k+1}^{j-1}\ts\bs_k\ts x_k^{i-j}
-x_k^{j-1}\ts\bs_k\ts x_{k+1}^{i-j}\ts\bigr)
\\
&=
\sum_{j=1}^i\ts\ts(-1)^{j-1}\ts
\bigl(\ts x_k^{j-1}\ts\bs_k\ts x_k^{i-j}
-x_k^{j-1}\ts\bs_k\ts x_k^{i-j}\ts\bigr)=0.
\endalign
$$
Now it follows directly from (\[2.2]) that for odd $i$ the sum
$x_1^i+\ldots+x_n^i$ commutes with $s_k$.
This sum then also commutes with $\bs_k$ due to (\[2.3])
and to Corollary 2.2\sq
\enddemo 

\np
It follows from the definition (\[1.0]) that for any
$b\in B(k,N)$ there is a unique element $b\pr\in B(k-1,N)$ such that
$$
\bs_k\ts b\ts\bs_k=b\pr\ts\bs_k\ts;
\Tag{2.43}
$$
cf. [W, Proposition 2.2]. Moreover, the map $b\mapsto b\pr$
evidently commutes with the left and right multiplication by elements
from the subalgebra $B(k-1,N)\subset B(k,N)$. In particular, due to
Lemma 2.1 we have 
$$
\bs_k\ts x_k^i\ts\bs_k=z_k^{(i)}\ts\bs_k\ts;\qquad i=0,1,2,\dots
\Tag{2.45}
$$
where $z_k^{(0)}=N$ and $z_k^{(1)},z_k^{(2)},\dots$ are 
central elements of the algebra $B(k-1,N)$. In Section 4 we will provide
explicit
formulas for these elements; see Corollary 4.3 and the subsequent remark. 
Here we will point out only some relations
that the definition (\[2.45]) implies.  

\proclaim{Lemma 2.5}
We have the relations
$$
-2\ts z_k^{(i)}=z_k^{(i-1)}+
\sum_{j=1}^i\ts\ts(-1)^j\ts z_k^{(i-j)}\ts z_k^{(j-1)};
\qquad
i=1,3,\dots\ts.
\Tag{2.475}
$$
\endproclaim

\demo{Proof}
Let us multiply the relation (\[2.425]) by $\bs_k$ on the left and on
the right. Then due to (\[1.1]) and (\[2.4]) by the definition (\[2.45])
we get
$$
\align
z_k^{(i)}\ts\bs_k
&=\bs_k\ts x_{k+1}^i\ts\bs_k
+
\sum_{j=1}^i\ts\ts
\bs_k\ts x_{k+1}^{j-1}\ts(\bs_k-1)\ts x_k^{i-j}\ts\bs_k
\Tag{2.4444}
\\
\qquad\ \ \ts\qquad
&
=(-1)^i\ts z_k^{(i)}\ts\bs_k
+
\sum_{j=1}^i\ts\ts
(-1)^{j-1}\ts\bigl(\ts z_k^{(j-1)}\ts z_k^{(i-j)}
-z_k^{(i-1)}\ts\bigr)\ts\bs_k\ts.
\endalign
\nopagebreak
$$
The last equality for odd $i$ implies (\[2.475])\sq
\enddemo

\np
From now on until the end of Section 3 we will assume that the parameter
$N$ is a positive integer. We will then have the chain of semisimple algebras
$$
C(1,N)\subset C(2,N)\subset\ldots\subset C(n,N).
\Tag{2.4775}
$$
Consider the subalgebra $X(n,N)$ in $C(n,N)$ generated by all the
central elements of $C(1,N),C(2,N),\dots,C(n,N)$. Each of the latter
algebras is semisimple.
So it follows from Corollary 1.4 that the subalgebra $X(n,N)$ is maximal
commutative.

There is a canonical basis in every representation space
$V(\la,n)$ of $C(n,N)$ corresponding to the chain (\[2.4775]); it
consists of the eigenvectors of the subalgebra $X(n,N)$. The basic
vectors are parametrized by the sequences  
$$
\La=\bigl(\La(1),\dots,\La(n)\bigr)\in\O(1,N)\times\ldots\times\O(n,N)
$$
where $\La(n)=\la$ and
each two neighbouring terms of the sequence differ by exactly one box.
Denote by $\L(\la,n)$ the set of all such sequences.
Let $v(\La)$ be the basic vector in $V(\la,n)$ corresponding to a
sequence $\La\in\L(\la,n)$.
Up to a scalar multiplier, it is uniquely determined by
the following condition:
$v(\La)\in V\bigl(\La(k),k\bigr)$ in the restriction
of $V(\la,n)$ onto $C(k,N)$ for any $k=1,\dots,n-1$.

We will regard $V(\la,n)$ as a representation of the algebra $B(n,N)$
also.
In the next section we will use the elements $x_1,\dots,x_n\in B(n,N)$
to describe the action of the generators
$s_1,\dots,s_{n-1};\bs_1,\dots,\bs_{n-1}$ of $B(n,N)$ on the vector
$v(\La)\in V(\la,n)$. It follows from Corollary 1.4 and Lemma 2.1 that
the images in $C(n,N)$ of the elements $x_1,\dots,x_n$ belong to the
subalgebra $X(n,N)$.
Denote by $x_k(\La)$ the eigenvalue of $x_k$ corresponding to the vector
$v(\La)$. For any $\La\in\L(\la,n)$ we will define $\La(0)$ as the
empty partition. 

\proclaim{Theorem 2.6}
Suppose that the diagrams $\La(k-1)$ and $\La(k)$ differ by the box
occuring in the row $i$ and the column $j$. Then
$$
x_k(\La)=\pm\left(\frac{N-1}2+j-i\right)
\Tag{2.5}
\nopagebreak
$$
where the upper sign in $\pm$ corresponds to the case
$\La(k)\supset\La(k-1)$ while the lower sign corresponds to
$\La(k)\subset\La(k-1)$.
\endproclaim

\demo{Proof}
If a box of the diagram $\la$ occurs in the row $i$ and the column $j$
then the difference $j-i$ is called the {\it content} of this box. Denote
by $n(\la)$ the number of the boxes in $\la$ and by $c(\la)$ the sum of
their contents. Due to Corollary 2.4 the element $x_1+\ldots+x_n$ is
central in $B(n,N)$. We shall prove that its eigenvalue in $V(\la,n)$ is
$$
c(\la,N)=\frac{N-1}2\ts n(\la)+c(\la).
$$
Applying this result to $k$ and $\La(k)$, $k-1$ and $\La(k-1)$ 
instead of $n$ and $\la$ we shall then obtain the equality (\[2.5]).

Consider $U^{\ot n}$ as a representation space of the algebra
$G\times B(n,N)$. Due to Proposition 1.2(b,c) we then have the
decomposition 
$$
U^{\ot n}=\underset{\la\in\O(n,N)}\to\oplus\ts U(\la,N)\ot V(\la,n)
\Tag{2.6}
$$
where $U(\la,N)$ is an irreducible representation of the group $G$.
In this decomposition $U(\la,N)$ does not depend on $n\geqslant n(\la)$;
see [Wy,\ts Theorem 5.7.F]. It suffices to demonstrate that for 
some vector $w\in U(\la,N)\ot V(\la,n)$
$$
(x_1+\ldots+x_n)\cdot w=c(\la,N)\ts w.
$$ 

Suppose that $n=n(\la)$. Due to [Wy, {\it loc.\ts cit.}]
any vector $w\in U(\la,N)\ot V(\la,n)$ is then {\it traceless:} we have
$\ol{(k,l)}\cdot w=0$ for all $k<l$. 
Thus $V(\la,n)$ is irreducible as a representatition of the group
$S(n)$; it is the representation corresponding to the diagram $\la$
[Wy,\ts Theorem 5.7.E]. Therefore
due to [Ma,\ts Examples I.1.3 and I.7.7]

$$
(x_1+\ldots+x_n)\cdot w=
\left(\ts\frac{N-1}2\ts n+
\sum_{1\leqslant k<l\leqslant n}\ts (k,l)\ts\right)\cdot w=
\left(\ts\frac{N-1}2\ts n+c(\la)\ts\right)\ts w
$$ 
so that the eigenvalue of $x_1+\ldots+x_n$ in $V(\la,n)$ for $n=n(\la)$
is equal to $c(\la,N)$. 

Now suppose that $n-n(\la)=2r>0$.
Then we will take $w=v_{\tts1}\ot v_{\tts2}^{\ot r}$ where
$$
v_{\tts1}\in U(\la,N)\ot V\bigl(\la,n(\la)\bigr)\subset U^{\ot n(\la)}
,\qquad
v_{\tts2}=\sum_{i=1}^N\ts u(i)\ot u(i)\ts\in U^{\ot2}.
$$
Let $E_{ij}$ with $i,j=1,\dots,N$ be the standard generators of the Lie
algebra $\frak{gl}(N)$. 
By definitions (\[1.5]) and (\[2.1])
the action in $U^{\ot n}$ of the element
$x_1+\ldots+x_n\in B(n,N)$ 
coincides with that of the Casimir element
$$
C=-\ts\frac14\ts\sum_{i,j=1}^N\ts(E_{ij}-E_{ji})^2
$$
of the universal enveloping algebra $\operatorname{U}(\frak{g})$
of the Lie algebra $\frak{g}=\frak{so}(N)\subset\frak{gl}(N)$. But 
$$
(E_{ij}-E_{ji})\cdot v_{\tts2}=0
$$
for any indices $i$ and $j\,$.
By the definition of the comultiplication on $\operatorname{U}(\frak{g})$
we now~get 
$$
\gather
(x_1+\ldots+x_n)\cdot w=C\cdot w=
(C\cdot v_{\tts1})\ot v_{\tts2}^{\ot r}=
\\
\bigl((x_1+\ldots+x_{n(\la)})\cdot v_{\tts1}\bigr)\ot v_{\tts2}^{\ot r}
=c(\la,N)\ts w\quad\square
\endgather
$$
\enddemo
 
\proclaim{Corollary 2.7}
Suppose that $N$ is odd or $N\geqslant2n-1$. Then\tts:
\itemitem{a)}
the images in $C(n,N)$
of the elements $x_1,\dots,x_n$ generate the algebra $X(n,N)$;
\itemitem{b)}
the images in $C(n,N)$ of the elements
$x_1^i+\ldots+x_n^i\ts$ with $i=1,3,\dots\ts$
generate the centre of the algebra $C(n,N)$.
\endproclaim

\demo{Proof}
Fix any diagram $\la\in\O(n,N)$. To prove the part a)
we have to demonstrate that for all different $\La\in\L(\la,n)$ the
 collections
of eigenvalues $\bigl(x_1(\La),\dots,x_n(\La)\bigr)$ are pairwise
distinct. Suppose that $\La,\Lap\in\L(n,N)$ are different, then
$\La(k-1)=\Lap(k-1)$ and $\La(k)\neq\Lap(k)$ for some $k\in\{2,\dots,n\}$.
Let $\La(k)$, $\Lap(k)$ differ from $\La(k-1)$ by the boxes occuring
in the rows $i,i\pr$ and the columns $j,j\pr$ respectively. 

If $\La(k),\Lap(k)\supset\La(k-1)$ or $\La(k),\Lap(k)\supset\La(k-1)$
then $j-i\neq j\pr-i\pr$ and $x_k(\La)\neq x_k(\Lap)$.
Suppose that $\La(k)\supset\La(k-1)$ and $\Lap(k)\subset\La(k-1)$, then
$$
x_k(\La)-x_k(\Lap)=N-1+j+j\pr-i-i\pr.      
$$
If $j+j\pr\geqslant3$ then $i+i\pr\leqslant N$ and
$x_k(\La)>x_k(\Lap)$. If $j\pr=j=1$ then $i\pr=i-1$ and
$x_k(\La)-x_k(\Lap)\neq0$ for the odd $N$ and for $N\geqslant2n-1>2i-2$.
This proves a).

Denote by $\C(\la)$ the unordered collection of the contents of all the
boxes of the diagram $\la$. Then $\la$ can be uniquely restored
from $\C(\la)$ since the boxes with the same content constitute the
diagonals of $\la$ and the lengths of the diagonals determine $\la$.
To prove the part b) we will show that the diagram $\la\in\O(n,N)$ can
be uniquely restored from the collection of the eigenvalues of the
elements $x_1^i+\ldots+x_n^i\ts$ with $i=1,3,\dots\ts$ in $V(\la,n)$.
Due to the equality of the formal power series in $u^{-1}$  
$$
\operatorname{exp}
\sum_{i=1,3,\dots}\ts2\ts(x_1^i+\ldots+x_n^i)\ts u^{-i}/i
\ts=\ts
\prod_{k=1}^n\ts\frac{u+x_k}{u-x_k}
$$
these eigenvalues determine the collection $\C(\la,N)$
obtained from $\C(\la)$ by removing all the pairs of the contents
$j-i,\ts j\pr-i\pr$ such that
$$
\frac{N-1}2+j-i=-\left(\frac{N-1}2+j\pr-i\pr\right).
$$
The latter condition implies that $j=j\pr=1$ and $i+i\pr=N+1$. Moreover,
then $i\neq i\pr$ and there is only one box in each of the rows $i,i\pr$
of the diagram $\la$. If $N\geqslant2n-1$ then $i+i\pr\leqslant2n-1<N+1$
so that $\C(\la)=\C(\la,N)$. 

Now let $N$ be odd. If the collection $\C(\la,N)$ does not
contain $(1-N)/2$ then there are less than $(N+1)/2$ rows in the diagram
$\la$ and $i+i\pr\leqslant N-2$ for any two different rows $i,i\pr$.
Then we have $\C(\la)=\C(\la,N)$ again. Suppose that $\C(\la,N)$ does
contain the number $(1-N)/2$. Then this number is minimal
in the collection $\C(\la,N)$ and occurs therein only once. Then we have
$$
\C(\la)=\C(\la,N)\sqcup\bigl\{\ts 1-i,i-N\ts\ts|\ts
\operatorname{min}
\left(
\C(\la,N)\setminus\bigl\{\frac{1-N}2\bigr\}
\right)
>1-i>
\frac{1-N}2\ts\bigr\}\,.
$$
Thus the collection $\C(\la)$ can be always restored from
$\C(\la,N)$. This proves b)\sq
\enddemo

\demo{Remark}
For $N=2,4,\dots,2n-2$
the parts a) and b) of Corollary 2.7 are not valid.
However, the elements $x_1,\dots,x_n$
will still suffice to describe the action in $V(\la,n)$ of the generators
$s_1,\dots,s_{n-1};\ts\bs_1,\dots,\bs_{n-1}$ of $B(n,N)$ 
for any positive integer $N$.
\enddemo

\heading
3. Young Orthogonal Form for $C(n,N)$
\endheading\reset

\np
It this section we will make explicit the matrix elements
$s_k(\La,\Lap),\bs_k(\La,\Lap)$ of the generators
$s_k,\bs_k\in B(n,N)$ in the canonical basis of the representation
$V(\la,n)$: 
$$
s_k\cdot v(\La)=\sum_{\Lap\in\L(\la,n)}\ts s_k(\La,\Lap)\ts v(\Lap),
\qquad
\bs_k\cdot v(\La)=
\sum_{\Lap\in\L(\la,n)}\ts\bs_k(\La,\Lap)\ts v(\Lap).
$$
Note that each of the vectors $v(\La)\in V(\la,n)$
here is defined up to a scalar
multiplier. Before specifying these multipliers we will determine
the diagonal matrix elements $s_k(\La,\La)$, $\bs_k(\La,\La)$ along
with all the products 
$s_k(\La,\Lap)\ts s_k(\Lap,\La)$,
$\bs_k(\La,\Lap)\ts\bs_k(\Lap,\La)$.

Let an index $k\in\{1,\dots,n-1\}$ and a sequence $\La\in\L(\la,n)$
be fixed.
Denote by $V(\La,k)$ the subspace in $V(\la,n)$ spanned by the
vectors $v(\Lap)$ such that $\Lap(l)=\La(l)$ for any $l\neq k$.
The action of $s_k$ and $\bs_k$ in $V(\la,n)$ preserves
this subspace.

\proclaim{Proposition 3.1}
Suppose that $\La(k-1)\neq\La(k+1)$. Then $\bs_k\cdot v(\La)=0$.
\endproclaim

\demo{Proof}
The diagrams $\La(k-1)$ and $\La(k+1)$ differ by two boxes; let
$j-i$ and $j\pr-i\pr$ be the contents of these boxes.
If $x_k(\La)+x_{k+1}(\La)\neq0$ then by applying to $v(\La)$
the first of the relations (\[2.4]) we obtain that
$\bs_k\cdot v(\La)=0$.

Now suppose that $x_k(\La)+x_{k+1}(\La)=0$.
Then by Theorem 2.6 we have $j=j\pr=1$ and $i+i\pr=N+1$. Therefore
either $\La(k-1)\subset\La(k+1)$ or $\La(k-1)\supset\La(k+1)$. In both
of these two cases the action of the elements $s_k$ and $\bs_k$ in
$V(\la,n)$ preserves the subspace $\CC\cdot v(\La)$. Moreover, then
we have $x_k(\La)=1/2$ and $x_{k+1}(\La)=-1/2$. Due to (\[1.1]) by
applying the first of the relations (\[2.3]) to the vector
$\bs_k\cdot v(\La)$ we obtain that $(N-2)\ts\bs_k\cdot v(\La)=0$.
This equality completes the proof for $N\neq2$.

Let $\ph$ and $\psi$ denote respectively the empty diagram
and the diagram consisting of two boxes in the first column. 
If $N=2$ then $\{\La(k-1),\La(k+1)\}=\{\ph,\psi\}$.
The representation $U(\psi,2)$ of 
$G=O(2,\CC)$ is the
determinant representation and
$$
U\bigl(\La(k+1),2\bigr)=
U\bigl(\La(k-1),2\bigr)\ot U(\psi,2).
$$
Therefore the action of $\bs_k$ in the space
$$
\align
\Hom_{B(k-1,2)}
&\bigl(V(\La(k-1),k-1)\ts,V(\La(k+1),k+1)\bigr)=
\\
\Hom_{\ts G\times B(k-1,2)}
&\bigl(U(\La(k+1),2)\ot V(\La(k-1),k-1)\ts,
U^{\ot(k+1)}\bigr)=
\\
\Hom_{\ts G}
&\bigl(U(\La(k+1),2)\ts,U(\La(k-1),2)\ot U^{\ot2}\bigr)
\endalign
$$
coincides with that of $\bs_1$ in $V(\psi,2)$.
This proves that $\bs_k\cdot v(\La)=0$ for $N=2$\sq
\enddemo

\proclaim{Proposition 3.2}
Let $\La(k-1)\neq\La(k+1)$. Then $x_k(\La)\neq x_{k+1}(\La)$ and
$$
s_k(\La,\La)=\bigl(x_{k+1}(\La)-x_k(\La)\bigr)^{-1}\ts.   
$$
\endproclaim

\demo{Proof}
By applying to the vector $v(\La)$ the second of the relations
(\[2.3]) we obtain that
$s_k(\La,\La)\ts\bigl(x_{k+1}(\La)-x_k(\La)\bigr)=1$\sq
\enddemo

\np
Observe that if $\La(k-1)\neq\La(k+1)$ then the space $V(\La,k)$ has
dimension at most two. Therefore due to the relation $s_k^2=1$
we get

\proclaim{Corollary 3.3}
Let $\La(k-1)\neq\La(k+1)$ and $v(\Lap)\in V(\La,k)$ with
$\La\neq\Lap$. Then
$$
s_k(\La,\Lap)\ts s_k(\Lap,\La)=
1-\bigl(x_{k+1}(\La)-x_k(\La)\bigr)^{-2}\ts.
$$
\endproclaim

\np
Two Young diagrams are {\it associated} if the sum of the
lengths of their first columns equals $N$ while the lengths of 
their other columns respecively coincide. In particular, for even $N$
a diagram is {\it self-associated} if its first column consists
of $N/2$ boxes. 

\proclaim{Lemma 3.4}
For any
$v(\Lap)\in V(\La,k)$
we have $x_k(\La)+x_{k}(\Lap)\neq0$
unless $N$ is odd and $\Lap=\La$ where the diagrams $\La(k-1),\La(k)$
are associated.
\endproclaim

\demo{Proof}
Let the diagrams $\La(k)$ and $\Lap(k)$ differ from $\La(k-1)$ by the
boxes with contents $j-i$ and $j\pr-i\pr$ respectively.
If either $\La(k)\subset\La(k-1)\subset\Lap(k)$ or
$\La(k)\supset\La(k-1)\supset\Lap(k)$ then $j-i\neq j\pr-i\pr$ and
$x_k(\La)+x_{k}(\Lap)\neq0$.

Assume now that either
$\La(k),\Lap(k)\subset\La(k-1)$ or $\La(k),\Lap(k)\supset\La(k-1)$.
Then the condition $x_k(\La)+x_{k}(\Lap)=0$ takes the form
$$
j-i+j\pr-i\pr=1-N
\Tag{3.01}
$$
which implies that $j=j\pr=1$. Then $i=i\pr$ due to our assumption.
Hence $\La=\Lap\,$. Moreover $N=2\tts i-1$
by (\[3.01]). So the diagrams $\La(k-1),\La(k)$ are associated
\sq
\enddemo

\newpage

\np
Let us now consider the case $\La(k-1)=\La(k+1)$. Due to Theorem 2.6
we then have $x_k(\Lap)+x_{k+1}(\Lap)=0$ for any $v(\Lap)\in V(\La,k)$.
The next two lemmas are contained in
[RW,\ts Theorem 2.4(b)]. We will include their short proofs here.

\proclaim{Lemma 3.5}
Suppose that $\La(k-1)=\La(k+1)$. Then
$$
\bs_k(\La,\La)=\frac
{\di\ts U\bigl(\La(k),N\bigr)}
{\di\ts U\bigl(\La(k+1),N\bigr)}\ts.
$$
\endproclaim

\demo{Proof}
Let $\tau_n$ denote the trace on the algebra $C(n,N)$ induced by the
usual matrix trace on $\operatorname{End}\bigl(U^{\ot n}\bigr)$.
We will also regard $\tau_n$ as a trace on the algebra $B(n,N)$.
Then due to the definition (\[1.0]) 
$$
\tau_{k+1}(\bs_k\ts b\ts\bs_k)=N\cdot\tau_k(b),\qquad b\in B(k,N).
$$
Let us apply this equality to an element $b\in B(k,N)$ such that for any
$\Lap\in\L(\la,n)$  
$$
b\cdot v(\Lap)=
\cases
v(\Lap)&\text{if }\quad\Lap(l)=\La(l)\ts\text{ for }l\leqslant k,
\\ 
0&\text{otherwise.}
\endcases
$$ 
Due to Proposition 3.1 we then obtain that
$$
N\cdot\di\ts U\bigl(\La(k+1),N\bigr)\ts\bs_k(\La,\La)=
N\cdot\di\ts U\bigl(\La(k),N\bigr)\ \ \square 
$$
\enddemo

\proclaim{Lemma 3.6}
Suppose that $\La(k-1)=\La(k+1)$. Then the image of the action of
$\bs_k$ in the subspace $V(\La,k)$ is one-dimensional.
\endproclaim

\demo{Proof}
Due to the relation $\bs_k^{\ts2}=N\ts\bs_k$ any eigenvalue of the
action of $\bs_k$ in the subspace $V(\La,k)$ is either $N$ or zero.
Lemma 3.5 along with the decomposition
$$
U\bigl(\La(k+1),N\bigr)\ot U=
\underset{v(\Lap)\in V(\La,k)}\to\oplus\ts U\bigl(\Lap(k),N\bigr)
$$
implies that the trace of this action is equal to $N$. Therefore the
eigenvalue $N$ of the action of $\bs_k$ in $V(\La,k)$ has
multiplicity one\sq  
\enddemo

\proclaim{Corollary 3.7}
Suppose that $\La(k-1)=\La(k+1)$ and $v(\Lap)\in V(\La,k)$. Then
$$
\bs_k(\La,\Lap)\ts\bs_k(\Lap,\La)=
\bs_k(\La,\La)\ts\bs_k(\Lap,\Lap)\ts.
$$
\endproclaim

\np
There are well known explicit formulas for the dimension of the
irreducible representation $U(\la,N)$ of the orthogonal group $G$;
see for instance [EK,\ts Section 3]. Due to Lemma 3.5
these formulas already provide certain expressions for the matrix element
$\bs_k(\La,\La)$. In this section we will employ the relations (\[2.3])
and (\[2.45]) to determine $\bs_k(\La,\La)$ independently
of any explicit formula for $\di\ts U(\la,N)$.

Suppose that $\La(k-1)=\La(k+1)=\mu$.
Let $l$ be the number of pairwise distinct rows (or columns)
in the diagram $\mu$. Then one can obtain $l+1$ diagrams
by adding a box to $\mu$ and $l$ diagrams by removing a box
from $\mu$.
Let $c_1,\dots,c_{l+1}$ and $d_1,\dots,d_l$ be the contents
of these boxes respectively.
Denote by $b_1,\dots,b_{2l+1}$
the numbers
$$
(N-1)/2+c_1\ts,\dots,(N-1)/2+c_{l+1}\ts,\ts
-(N-1)/2-d_1\ts,\dots,-(N-1)/2-d_l
$$
taken in an arbitrary order. Introduce the formal power series in $u^{-1}$
$$
Q(\mu,u)=
\sum_{i\geqslant0}\ts q_i(\mu)\ts u^{-i}=
\ts\prod_{j=1}^{2l+1}\ts\ts\frac{u+b_j}{u-b_j}\ts;
\Tag{3.*}
$$
the coefficients $q_1(\mu),q_2(\mu),\dots$
are the symmetric {\it Schur $q$-functions} in $b_1,\dots,b_{2l+1}$.

Further, let $m$ be the total number of boxes in the diagram $\mu$.
Let $e_1,\dots,e_{m}$ be the contents of all these boxes. 
Denote by $a_1,\dots,a_m$
the numbers
$$
(N-1)/2+e_1\ts,\dots,(N-1)/2+e_m
$$
taken in an arbitrary order. 
Then we have another expression for the series (\[3.*]).

\proclaim{Lemma 3.8}
We have the equality
$$
Q(\mu,u)
=
\frac{u+(N-1)/2}{u-(N-1)/2}
\ts\cdot\ts\ts
\prod_{j=1}^{m}\ts\ts
\frac{(u+a_j)^2-1}{(u-a_j)^2-1}
\ts
\frac{(u-a_j)^2}{(u+a_j)^2}\ts\ts.
$$
\endproclaim

\demo{Proof}
For any $h\in\CC$ and $i\geqslant0$ we have the equality
$$
\sum_{j=1}^{l+1}\ts(h+c_j)^i\ts-\ts\sum_{j=1}^l\ts(h+d_j)^i=
h^i+\sum_{j=1}^m\ts(h+e_j+1)^i-2(h+e_j)^i+(h+e_j-1)^i,
$$
it can be verified by induction on $m$.
Let us multiply each side of this equality by $2u^{-i}/i$ 
and take the sum over all odd~$i$. Exponentiate
the resulting sums. When $h=(N-1)/2$ we get the required statement\sq
\enddemo

\np
Denote by $z_k^{(i)}(\mu)$ the eigenvalue of the central element
$z_k^{(i)}\in B(k-1,N)$ defined by (\[2.45]) in the representation
$V(\mu,k-1)$. Introduce the generating function
$$
Z(\mu,u)=\sum_{i\geqslant0}\ts z_k^{i}(\mu)\ts u^{-i}.
$$

\proclaim{Theorem 3.9}
We have the equality
$$
Z(\mu,u)=\bigl(u+1/2\bigr)\cdot Q(\mu,u)-u+1/2\ts.
$$
\endproclaim

\demo{Proof}
Consider the generating series
$$
Z_k(u)=\sum_{i\geqslant0}\ts z_k^{(i)}\ts u^{-i}\ts\in\ts B(k-1,N)[[u^{-1}]].
\Tag{2.455555}
$$
Then determine the series $Q_k(u)\in B(k-1,N)[[u^{-1}]]$ by the equality
$$
Z_k(u)=\bigl(u+1/2\bigr)\cdot Q_k(u)-u+1/2\ts.
\Tag{2.4666666}
$$
We have to prove that the eigenvalue of $Q_k(u)$ in 
$V(\mu,k-1)$ is exactly $Q(\mu,u)$. Later on in a more general
setting we shall prove the equality
$$
Q_k(u)
=
\frac{u+(N-1)/2}{u-(N-1)/2}
\ts\cdot\ts\ts
\prod_{l=1}^{k-1}\ts\ts
\frac{(u+x_l)^2-1}{(u-x_l)^2-1}
\ts
\frac{(u-x_l)^2}{(u+x_l)^2}\ts\ts;
\Tag{3.**}
$$
see Corollary 4.3. Due to (\[3.**])
the required statement follows from
Theorem~2.6 and Lemma 3.8\sq 
\enddemo

\proclaim{Corollary 3.10}
Suppose that $\La(k-1)=\La(k+1)=\mu$ and let $x_k(\La)=b\ts$. Then
$$
\bs_k(\La,\La)=
\left\{\aligned
(\ts2b+1\ts)\ts\prod\Sb b_j\neq b\endSb
\ts\frac{b+b_j}{b-b_j}
&\quad\text{if}\quad b\neq-1/2\ts;
\\
-\ \prod\Sb b_j\neq b\endSb
\ts\frac{b+b_j}{b-b_j}
&\quad\text{if}\quad b=-1/2\ts.
\endaligned\right.
$$
\endproclaim

\demo{Proof}
Let us make use of the relation
$$
\bs_k\cdot Z_k(u)\ts u^{-1}=\bs_k\ts(u-x_k)^{-1}\ts\bs_k\ts.
\Tag{3.15}
$$
Any eigenvalue of $x_k$ in $V(\La,k)$ distinct from $-1/2$ has
multiplicity one. Moreover, it then appears only once in
the collection $b_1,\dots,b_{2l+1}$. Therefore if $b\neq-1/2$~we~get
from (\[3.15])
$$
\bs_k(\La,\La)=\underset{u=b}\to{\operatorname{res}}\ts\ts
{Z(\mu,u)}/{u}=(\ts2b+1\ts)\ts\prod_{b_j\neq b}\ts
\frac{b+b_j}{b-b_j}\ts;    
$$
here we use Theorem 3.9 and the inequality
$\bs_k(\La,\La)\neq0$ provided by Lemma 3.5.  

If $b=-1/2$ appears as an eigenvalue of $x_k$ in $V(\La,k)$ it
has multiplicity two. Then it appears twice in the
collection $b_1,\dots,b_{2l+1}$.
Let $x_k(\La)=x_k(\La\pr)=-1/2$ where
$v(\Lap)\in V(\La,k)$ and $\La(k)\neq\Lap(k)$. Then the diagrams
$\La(k)$ and $\Lap(k)$ are associated. Representations
$U(\La(k),N)$ and $U(\Lap(k),N)$ of $G$ 
have the same dimension [Wy\ts,\ts Theorem 5.9.A]. So 
$\bs_k(\La,\La)=\bs_k(\Lap,\Lap)\neq0$ by Lemma~3.5. Therefore
$$
\bs_k(\La,\La)=\underset{u=b}\to{\operatorname{res}}\ts\ts
{Z(\mu,u)}/{2u}=-\ts\prod_{b_j\neq b}\ts
\frac{b+b_j}{b-b_j}
\quad\square    
$$
\enddemo

\proclaim{Proposition 3.11}
Suppose that $\La(k-1)=\La(k+1)$ and $v(\Lap)\in V(\La,k)$. Then
$$
s_k(\La,\Lap)=
\bigl(\bs_k(\La,\Lap)-\de(\La,\Lap)\bigr)\ts
\bigl(x_k(\La)+x_k(\Lap)\bigr)^{-1}
\Tag{3.3}
$$
unless $N$ is odd and $\Lap=\La$ where the diagrams $\La(k),\La(k-1)$
are associated. In the latter case $s_k(\La,\La)=1$.
\endproclaim

\demo{Proof}
By the equality $x_k(\Lap)+x_{k+1}(\Lap)=0$ the first of the
relations (\[2.3]) implies
$$
s_k(\La,\Lap)\ts(x_k(\La)+x_k(\Lap)\bigr)=
\bigl(\bs_k(\La,\Lap)-\de(\La,\Lap)\bigr).
$$
Thus when $x_k(\La)+x_k(\Lap)\neq0$ we obtain the equality (\[3.3]).

\smallskip

Now assume that $x_k(\La)+x_k(\Lap)=0$. Then by Lemma 3.4 the number
$N$ is odd and $\Lap=\La$ where the diagrams $\La(k),\La(k-1)$
are associated. Then $\bs_k(\La,\La)=1$ by Lemma 3.5 and
$x_k(\La)=0$. Moreover, all the eigenvalues of the element $x_k$
in $V(\La,k)$ are then pairwise distinct.
Consider the diagonal matrix element of the
relation $s_k\ts\bs_k=\bs_k$ in $V(\La,k)$ corresponding
to the vector $v(\La)$. By making use of the equality (\[3.3]) and of
Corollary 3.7 we obtain that
$$
s_k(\La,\La)+\sum_{\Lapp\neq\La}
\ts\bs_k(\Lapp,\Lapp)/x_k(\Lapp)=1
\Tag{3.4}
$$
where $v(\Lapp)\in V(\La,k)$.
If $x_k(\Lapp)=b\neq0$ then
$$
\bs_k(\Lapp,\Lapp)/x_k(\Lapp)=
\underset{u=b}\to{\operatorname{res}}\ts\ts
{Z(\mu,u)}/{u^2}\ts
$$
while by Theorem 3.9 the residues of ${Z(\mu,u)}/{u^2}$
at $u=0,\infty$ equal zero.
Therefore
$$
\sum_{\Lapp\neq\La}\ts\bs_k(\Lapp,\Lapp)/x_k(\Lapp)=0.
$$
By comparing this equality with (\[3.4]) we complete the proof of
Proposition 3.11
\sq
\enddemo 

\np
Now let the index $k$ run through the set $\{1,\dots,n-1\}$ while
the sequences $\La,\Lap$ run through the set $\L(\la,n)$.
If $v(\Lap)\notin V(\La,k)$ then $s_k(\La,\Lap)=\bs_k(\La,\Lap)=0$.

Suppose that $v(\Lap)\in V(\La,k)$.
As we have already mentioned, the vectors $v(\La),v(\Lap)\in V(\la,n)$
are defined up to scalar multipliers. Up to the choice of these
multipliers Proposition 3.1 and Corollaries 3.7,\ts3.10 describe the matrix
element $\bs_k(\La,\Lap)$ while Propositions 3.2,\ts3.10 and Corollary 3.3
describe the matrix element $s_k(\La,\Lap)$.
The following theorem
completes the description of these matrix elements.

\proclaim{Theorem 3.12}
Suppose that $v(\Lap)\in V(\La,k)$ and $\La\neq\Lap$. Then one can
assume:
$$
\align
s_k(\La,\Lap)=s_k(\Lap,\La)>0\quad&\text{if}\quad\La(k-1)\neq\La(k+1),
\Tag{3.33}
\\
\bs_k(\La,\Lap)=\bs_k(\Lap,\La)>0\quad&\text{if}\quad\La(k-1)=\La(k+1).
\Tag{3.66}
\endalign
$$
\endproclaim

\demo{Proof}
Let us demonstrate first that for all $k$ and $\La,\Lap$
we can assume the equalities
$$
s_k(\La,\Lap)=s_k(\Lap,\La)
\quad\text{and}\quad
\bs_k(\La,\Lap)=\bs_k(\Lap,\La).
\Tag{3.5}
$$
Consider the non-degenerate symmetric bilinear form $\fm$
on $U^{\ot n}$ which is the product
of the standard forms on the factors $U$. The action of the group $G$ in
$U^{\ot n}$ preserves $\fm$ and the direct summands in (\[2.6])
are othogonal with respect to $\fm$. The restriction of $\fm$
onto the direct summand $U(\la,N)\ot V(\la,n)$ splits into the product
of a $G$\ts-invariant bilinear symmetric form on $U(\la,N)$ and of 
a certain form 
on $V(\la,N)$. The vectors $v(\La)\in V(\la,n)$ are
orthogonal with respect to the latter form. The
actions of $s_k$ and $\bs_k$ in $U^{\ot n}$ are
self-adjoint with respect to $\fm$. 
Therefore the vectors $v(\La)$ can be so chosen that the equalities
(\[3.5]) hold.

We will now assume that the equalities in (\[3.33]) and (\[3.66])
do hold. The matrix elements in these equalities then belong to
$\RR\setminus\{0\}$ by Corollaries 3.3 and 3.7.  We have to prove
that the vectors $v(\La)$ can be so chosen that the inequalities in
(\[3.33]) and (\[3.66]) also hold. 
Let an index $k\in\{1,\dots,n-2\}$ and a sequence $\La\in\L(\la,n)$
be fixed.
Denote by $V(\La,k,k+1)$ the subspace in $V(\la,n)$ spanned by the
vectors $v(\Lap)$ such that $\Lap(l)=\La(l)$ for any $l\neq k,k+1$.
The action of $s_k\ts, s_{k+1}$ and $\bs_k\ts,\bs_{k+1}$ in
$V(\la,n)$ preserves this subspace.

Note that due to Proposition 1.1 the generators $s_k\ts,\bs_k$
of the algebra $B(n,N)$ are {\it local\ts} in the sense [V1,V2]:
the only relations between $s_k\ts,\bs_k$ and $s_l\ts,\bs_l$ with
$|k-l|>1$ are the commutation relations (\[1.4]). 
Therefore it suffices to choose only the vectors $v(\Lap)\in V(\La,k,k+1)$
so that for every two distinct $v(\Lap),v(\Lapp)\in V(\La,k,k+1)$
$$
\align
s_k(\Lap,\Lapp)>0\quad
\text{if}\quad
&\Lap(k+1)=\Lapp(k+1)\neq\La(k-1),
\Tag{3.501}
\\
s_{k+1}(\Lap,\Lapp)>0\quad
\text{if}\quad
&\Lap(k)=\Lapp(k)\neq\La(k+2),
\Tag{3.502}
\\
\bs_k(\Lap,\Lapp)>0\quad
\text{if}\quad
&\Lap(k+1)=\Lapp(k+1)=\La(k-1),
\Tag{3.503}
\\
\bs_{k+1}(\Lap,\Lapp)>0\quad
\text{if}\quad
&\Lap(k)=\Lapp(k)=\La(k+2).
\Tag{3.504}
\endalign
$$

The diagrams $\La(k-1)$ and $\La(k+2)$ may differ either by one or by
three boxes. Let us consider the latter case; cf.  [Mo].
In this case we have
$\Lap(k+1)\neq\La(k-1)$ and $\Lap(k)\neq\La(k+2)$ for any vector
$v(\Lap)\in V(\La,k,k+1)$. Furthermore, then
$\di\ts V(\La,k,k+1)\in\{1,3,6\}$. We will treat each of these
possibilities separately.
\smallskip
\np
i) If $\di\ts V(\La,k,k+1)=1$ then we can choose the vector $v(\La)$
arbitrarily.
\smallskip
\np
ii) Suppose that $\di\ts V(\La,k,k+1)=3$. Initially let us make an
arbitrary choice of the
basic vectors $v(\La),v(\Lap),v(\Lapp)\in V(\La,k,k+1)$.
We can assume that
$$
\La(k)\neq\Lap(k)=\Lapp(k)
\quad\text{and}\quad
\Lap(k+1)\neq\La(k+1)=\Lapp(k+1).
$$
This assumption is illustrated by the following graph where each edge 
indicates a difference by exactly one box\hskip.5pt:
$$
\xymatrix{
\,\La(k-1)\,\ar@{-}[r]\ar@{-}[dr]
&
\,\La(k)\,\ar@{-}[r]
&
\,\La(k+1)\,\ar@{-}[r]
&
\,\La(k+2)\,
\\
&
\,\Lap(k)\,\ar@{-}[r]\ar@{-}[ur]\ar@{-}
&
\,\Lap(k+1)\,\ar@{-}[ur]
&
}
$$
Then it suffices to take instead of 
$v(\Lap)$ and $v(\Lapp)$ respectively the vectors 
$$
v(\Lap)\cdot\si\ts s_k(\La,\Lapp)\ts s_{k+1}(\Lap,\Lapp)
\quad\text{and}\quad
v(\Lapp)\cdot\si\ts s_k(\La,\Lapp).
$$
\np
iii) Suppose that $\di\ts V(\La,k,k+1)=6$. Initially let us again
make an arbitrary choice of the
basic vectors in $V(\La,k,k+1)$
$$
v(\La),
\ 
v(\Om),
\ 
v(\Omp),
\ 
v(\Lap),
\ 
v(\Lapp),
\ 
v(\Ompp).
\Tag{3.6}
$$
We can assume that
$$
\alignat2
\Om(k)&=\Lap(k),
&
\Om(k+1)&=\La(k+1),
\\
\Omp(k)&=\La(k),
&
\Omp(k+1)&=\Lapp(k+1),
\\
\Ompp(k)&=\Lapp(k),\qquad
&
\Ompp(k+1)&=\Lap(k+1).
\endalignat
$$
This assumption is illustrated by the next graph\hskip.5pt:
$$
\xymatrix{
&
\,\La(k)\,\ar@{-}[dr]\ar@{-}[r]
&
\,\La(k+1)\,\ar@{-}[dr]
&
\\
\,\La(k-1)\,\ar@{-}[r]\ar@{-}[dr]\ar@{-}[ur]
&
\,\Lap(k)\,\ar@{-}[ur]\ar@{-}[dr]
&
\,\Lapp(k+1)\,\ar@{-}[r]
&
\,\La(k+2)\,
\\
&
\,\Lapp(k)\,\ar@{-}[ur]\ar@{-}[r]
&
\,\Lap(k+1)\,\ar@{-}[ur]
&
}
$$
By applying the first
relation in (\[1.2]) to the vector $v(\La)$ and taking the coefficient
at $v(\Ompp)$ we get the equality
$$
s_k(\La,\Om)\ts s_{k+1}(\Lap,\Om)\ts s_k(\Lap,\Ompp)=
s_{k+1}(\La,\Omp)\ts s_k(\Lapp,\Omp)\ts s_{k+1}(\Lapp,\Ompp).
$$
Now it suffices to take instead of the vectors (\[3.6])
respectively the vectors 
$$
\gather
v(\La),
\quad 
v(\Om)\cdot\si\ts s_k(\La,\Om),
\quad
v(\Omp)\cdot\si\ts s_{k+1}(\La,\Omp),
\\
v(\Lap)\cdot\si\ts s_k(\La,\Om)\ts s_{k+1}(\Lap,\Om),
\quad
v(\Lapp)\cdot\si\ts s_{k+1}(\La,\Omp)\ts s_k(\Lapp,\Omp)
\\
v(\Ompp)\cdot\si\ts s_k(\La,\Om)\ts s_{k+1}(\Lap,\Om)\ts s_k(\Lap,\Ompp).
\endgather
$$

Finally, let us consider the case when the diagrams 
$\La(k-1)$ and $\La(k+2)$ differ by only one box.
Since $V(\La,k,k+1)=V(\Lap,k,k+1)$ for any $\Lap\in\L(\la,n)$
such that $\La(l)=\Lap(l)$ with $l\neq k,k+1$ we can assume that
$$
\La(k)=\La(k+2)
\quad\text{and}\quad
\La(k+1)=\La(k-1).
\Tag{3.7}
$$
Let us make an arbitrary initial choice of every basic vector
$v(\Lap)\in V(\La,k,k+1)$. 
Consider any vector $v(\Om)\in V(\La,k,k+1)$ such that
$$
\Om(k)\neq\La(k)
\quad\text{and}\quad
\Om(k+1)\neq\La(k+1).
$$
Then  consider the vectors 
$v(\Omp),v(\Ompp)\in V(\La,k,k+1)$ such that
$$
\alignat2
\Omp(k)&=\La(k),
&
\Omp(k+1)&=\Om(k+1),
\\
\Ompp(k)&=\Om(k),\qquad
&
\Ompp(k+1)&=\La(k+1).
\endalignat
$$
This assumption is illustrated by another graph\hskip.5pt: 
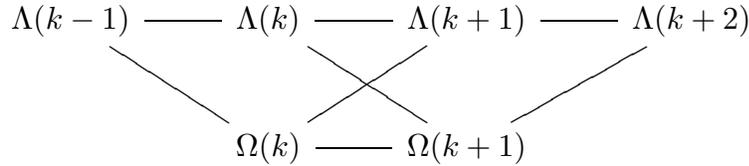
$$
\xymatrix{
\,\La(k-1)\,\ar@{-}[r]\ar@{-}[dr]
&
\,\La(k)\,\ar@{-}[r]\ar@{-}[dr]
&
\,\La(k+1)\,\ar@{-}[r]
&
\,\La(k+2)\,
\\
&
\,\Om(k)\,\ar@{-}[r]\ar@{-}[ur]\ar@{-}
&
\,\Om(k+1)\,\ar@{-}[ur]
&
}
$$
By applying the first relation in (\[1.3]) to the vector $v(\La)$
and taking the coefficient at $v(\Om)$ we get the equality
$$
\bs_k(\La,\La)\ts\bs_{k+1}(\La,\Omp)\ts s_k(\Om,\Omp)=
\bs_k(\La,\Ompp)\ts s_{k+1}(\Om,\Ompp).\hskip36pt
$$
Since $\bs_k(\La,\La)>0$ by Lemma 3.5, this equality implies that
$$
\si\ts\bs_{k+1}(\La,\Omp)\ts s_k(\Om,\Omp)=
\si\ts\bs_k(\La,\Ompp)\ts s_{k+1}(\Om,\Ompp).
$$
We will keep to the initial choice of the vector $v(\La)$
and replace each $v(\Om)$ by
$$
v(\Om)\cdot\si\bs_{k+1}(\La,\Omp)\ts s_k(\Om,\Omp)
\Tag{3.8}
$$
For $\Lap\neq\La$ where
$\Lap(k)=\La(k)$ or $\Lap(k+1)=\La(k+1)$
we will replace $v(\Lap)$ by
$$
\alignat3
v(\Lap)&\cdot\si\ts\bs_k(\La,\Lap)
&
\qquad&\text{if}\quad
&
\Lap(k+1)&=\La(k+1),
\\
v(\Lap)&\cdot\si\ts\bs_{k+1}(\La,\Lap)
&
\qquad&
\text{if}\quad
&
\Lap(k)&=\La(k).
\endalignat
$$
Due to Corollary 3.7 the latter replacement will make
all the matrix elements
of $\bs_k$ in $V(\La,k)$ and $\bs_{k+1}$ in $V(\La,k+1)$ positive. 
But for any $v(\Lap)\in V(\La,k,k+1)$
$$
\align
\Lap(k+1)=\La(k-1)
&\quad\Rightarrow\quad
v(\Lap)\in V(\La,k),
\\
\Lap(k)=\La(k+2)
&\quad\Rightarrow\quad
v(\Lap)\in V(\La,k+1)
\endalign
$$
due to (\[3.7]). So the inequalities (\[3.503]),\ts (\[3.504]) in 
$V(\La,k,k+1)$ will be then satisfied.

\np
Furthermore, for any two distinct $v(\Lap),v(\Lapp)\in V(\La,k,k+1)$
$$
\gather
\Lap(k+1)=\Lapp(k+1)\neq\La(k-1)\quad\Rightarrow\quad
\Lap(k)=\La(k)\ \ \text{or}\ \ \Lapp(k)=\La(k),
\\
\Lap(k)=\Lapp(k)\neq\La(k+2)\quad\Rightarrow\quad
\Lap(k+1)=\La(k+1)\ \ \text{or}\ \ \Lapp(k+1)=\La(k+1).
\endgather
$$
So the replacement of each $v(\Om)$ by the respective
vector (\[3.8]) will make the inequalities (\[3.501]),\ts (\[3.502])  
in $V(\La,k,k+1)$ satisfied. Theorem 3.12 is proved\sq
\enddemo

\heading
4. Affine Brauer Algebra
\endheading\reset

\np
In this section again we assume that $N$ is an arbitrary
complex number. We will now use the results of Section 2 as
a motivation to introduce a new object. This is the complex
associative algebra generated by the algebra $B(n,N)$ along with
the pairwise commuting elements $y_1,\dots,y_n$ 
and central elements $w_1,w_2,\ts\dots$
subjected to the following relations. We impose the relations
$$
\align
s_k\ts y_l=y_l\ts s_k\ts,
&\quad
\bs_k\ts y_l=y_l\ts\bs_k\ts;
\qquad
l\neq k,k+1;
\Tag{4.2}
\\
s_k\ts y_k-y_{k+1}\ts s_k=\bs_k-1\ts,
&\quad
s_k\ts y_{k+1}-y_k\ts s_k=1-\bs_k\ts;
\Tag{4.3}
\\
\bs_k\ts(y_k+y_{k+1})=0\ts,
&\quad
(y_k+y_{k+1})\ts\bs_k=0\ts.
\Tag{4.4}
\\
\intertext{Moreover, we impose the relations}
\bs_1\ts y_1^i\ts\bs_1=w_i\ts\bs_1\ts;
&\qquad
i=1,2,\ts\ldots\ts.
\Tag{4.45}
\endalign
$$
We view this algebra as an analogue of the degenerate affine Hecke
algebra $\Hen$ considered in [C1,C2] and [D]; see Corollary 4.9 below.
We will denote this algebra by $\Wen$ and call it the
{\it  affine Brauer algebra} here. Initially we proposed to call it the {\it degenerate affine Wenzl algebra} in  honour of H.\ts Wenzl
who used the maps
$$
B(k,N)\to B(k-1,N):\ts\ts b\ts\mapsto\ts b^\prime\ts;
\qquad
k=1,2,\dots,n-1
$$
defined by (\[2.43]) to prove that
the algebra $B(n,N)$ is semisimple when $N$ is not an integer [W].
However, this initial terminology has not been successful.

It is convenient to put $w_{0}=N$.
The equality (\[4.45]) is then valid for $i=0$ also.
The assignments
$$
y_k\mapsto x_k,
\quad
w_i\mapsto z_1^{(i)}
$$
define a homomorphism
$$
\pi:\Wen\to B(n,N)
\Tag{4.455555}
$$
identical on $B(n,N)$, see
the relations (\[2.2]),(\[2.3]),(\[2.4]),(\[2.45]).
Note that
in the proofs of Corollary 2.4 and Lemma 2.5 we used not the definition
(\[2.1]) of the elements $x_1,\dots,x_n$ but the latter relations.
Therefore the relations (\[4.2]) to (\[4.45]) imply that
$$
-2\ts w_i=w_{i-1}+
\sum_{j=1}^i\ts\ts(-1)^j\ts w_{i-j}\ts w_{j-1};
\qquad
i=1,3,\ts\ldots.
\Tag{4.475}
$$
In particular, we have $w_{1}=N(N-1)/2$. Morover, the following
proposition holds.

\proclaim{Proposition 4.1}
The elements $y_1^i+\ldots+y_n^i$ with $i=1,3,\ldots$
are central in the algebra $\Wen$.
\endproclaim 

\np
Due to the defining relations (\[4.2]) to (\[4.45])
we have an ascending chain of sublgebras
$$
A\hskip.5pt(1,N)
\subset
A\hskip.5pt(2,N)
\subset
\ldots
\,.
$$

\proclaim{Proposition 4.2}
For each $k=1,2,\ts\ldots$
we have the equalities
$$
\bs_k\ts y_k^i\ts\bs_k=w^{(i)}_k\ts\bs_k\ts;
\qquad
i=0,1,2,\ts\ldots
\Tag{4.5}
$$
where $w^{(i)}_k$ is a central element of
the algebra $A\hskip.5pt(k-1,N)$.
The generating series
$$
W_k(u)=\sum_{i\geqslant0}\ts w_k^{(i)}\ts u^{-i}
$$
\vskip-12pt
\np satisfy
\vskip-6pt
$$
\qquad
\frac
{W_{k+1}(u)+u-1/2}
{\hskip9pt W_k(u)+u-1/2}
=
\frac{(u+y_k)^2-1}{(u-y_k)^2-1}
\cdot
\frac{(u-y_k)^2}{(u+y_k)^2}\ts\ts.
\Tag{4.6}
$$
\endproclaim

\demo{Proof}
We use the induction on $k$.
The equalities (\[4.5])
hold for $k=1$ and $w^{(i)}_1=w_{i}$ by definition.
Assume that the equalities (\[4.5]) are valid for $k=1,\dots,n$
and that the corresponding series $W_1(u),\dots,W_n(u)$ satisfy (\[4.6]).
Due to Proposition 4.1 it then suffices to verify
the equalities (\[4.5]) and (\[4.6]) for $k=n+1$ and $k=n$ respectively.
We will work with the formal power series in $u^{-1}$
$$
\sum_{i\geqslant0}\ts y_k^{i}\ts u^{-i}
=\frac{u}{u-y_k}\ts;
\qquad
k=n,n+1\ts.
$$
We will also use the following corollary to the first relation in (\[4.3])
and the second relation in (\[4.4])
with $k=n\ts$:
$$
s_n\ts
\frac1{u-y_n}=
\frac1{u-y_{n+1}}\ts s_n+
\frac1{u+y_n}\ts\bs_n\frac1{u-y_n}-
\frac1{(u-y_n)\ts(u-y_{n+1})}\ts.
\Tag{4.7}
$$
By multiplyng the equality (\[4.7]) on the left by $\bs_n$ and replacing
$u$ by $-u$ we get
$$
\bs_n\ts\frac1{u-y_n}\ts s_n=
\bs_n\ts\frac1{u+y_n}+
\bs_n\ts\ts\frac{W_n(u)}{u(u+y_n)}-
\bs_n\ts\frac1{u^2-y_n^2}
\Tag{4.8}
$$
by the inductive assumption.
Let us now multiply the equality (\[4.7]) on the right by $s_n$ and use
(\[4.7]) once more along with (\[4.8]). We then obtain the equalities
$$
\align
s_n\ts\frac1{u-y_n}\ts s_n\ts
=&\ts
\frac{(u-y_n)^2-1}{(u-y_n)^2\ts(u-y_{n+1})}-
\frac1{u-y_n}\ts s_n\ts\frac1{u-y_n}+
\frac1{u+y_n}\ts\bs_n\ts\frac1{u+y_n}\ts+
\\
&\ts
\frac1{u^2-y_n^2}\ts\bs_n\ts\frac1{u-y_n}+
\frac1{u+y_n}\ts\bs_n\ts\ts\frac{W_n(u)}{u(u+y_n)}-
\frac1{u+y_n}\ts\bs_n\ts\frac1{u^2-y_n^2}
\ts.
\endalign
$$
By multiplying the last part of these equalities by $\bs_{n+1}$ on the left
and right we get 
$$
\ts
\frac{(u-y_n)^2-1}{(u-y_n)^2}
\cdot
\bs_{n+1}\ts\frac1{u-y_{n+1}}\ts\bs_{n+1}
+
\frac{W_n(u)}{u(u+y_n)^2}\ts\bs_{n+1}
+
\frac{2\ts(1-2u)\ts y_n}{(u^2-y_n^2\ts)^2}\ts\bs_{n+1}
\Tag{4.9}
$$
due to the relations
$
\bs_{n+1}\ts y_n=y_n\ts\bs_{n+1}
$
and
$
\bs_{n+1}\ts s_n\ts\bs_{n+1}=\bs_{n+1}\ts\bs_n\ts\bs_{n+1}=\bs_{n+1}\ts.
$
On the other hand, since 
$
s_{n+1}\ts y_n\ts s_{n+1}=y_n
$
we have by the inductive assumption 
$$
\bs_{n+1}\ts s_n\cdot\frac1{u-y_n}\cdot s_n\ts\bs_{n+1}=
\frac{W_n(u)}u\cdot\bs_{n+1}
\ts.
$$
By comparing the the last part of these equalities with the expression
(\[4.9]) we get
$$ 
\align
&\frac{(u-y_n)^2-1}{(u-y_n)^2}
\cdot
\bs_{n+1}\cdot\frac{u}{u-y_{n+1}}\cdot\bs_{n+1}
=
\\
&\frac{(u+y_n)^2-1}{(u+y_n)^2}
\cdot
W_n(u)\ts\bs_{n+1}-
\frac{2u\ts(1-2u)\,y_n}{(u^2-y_n^2\ts)^2}\cdot\bs_{n+1}
\ts.
\endalign
$$
The latter equality shows that
$$
\bs_{n+1}\cdot\frac{u}{u-y_{n+1}}\cdot\bs_{n+1}=
W_{n+1}(u)\cdot\bs_{n+1}
$$
where the series $W_{n+1}(u)$ satisfies (\[4.6]) for $k=n$.
Proposition 4.2 is proved
\ $\square$
\enddemo 

\np
Consider the series $Z_k(u)$ and $Q_k(u)$ defined by
(\[2.455555]) and (\[2.4666666]) respectively. Since $x_1=(N-1)/2$
we have
$$
Q_1(u)=\frac{u+(N-1)/2}{u-(N-1)/2}\ts.
$$
Furthermore,
$$
\pi:W_k(u)\mapsto Z_k(u)
$$
for every $k=1,2,\ts\ldots$ by (\[4.5]). Thus we obtain a
corollary to Proposition 4.2.

\proclaim{Corollary 4.3}
For every $k=1,2,\ts\ldots$ the equality {\rm(\[3.**])} holds.
\endproclaim

\np
In the remaining part of this section will construct a linear basis
in the algebra $\Wen$.
Let us equip the algebra $\Wen$ with an ascending filtration
by defining the degrees of its generators in the following way:
$$
\deg s_k=\deg\bs_k=0,
\quad
\deg y_k=1,
\quad
\deg w^{(i)}=0.
$$
Denote by $\ty_k$ the image of the element $y_k\in\Wen$
in the corresponding graded algebra $\gr\Wen$. 
In the latter algebra by the relations (\[4.2]) to (\[4.4]) we have
$$
\align
s\ts\ty_k\ts s^{-1}=\ty_{s(k)}\ts,
&\quad
s\in S(n).
\Tag{4.12}
\\
\intertext{These relations along with (\[4.2]) and (\[4.4]) imply that}
\ol{(k,l)}\ts\ts\ty_m=\ty_m\ts\ts\ol{(k,l)}\ts,
&\quad
m\neq k,l\ts;
\Tag{4.13}
\\
\ol{(k,l)}\cdot(\ty_k+\ty_l)=0\ts,
&\quad
(\ty_k+\ty_l)\cdot\ol{(k,l)}=0\ts;
\qquad
k\neq l.
\Tag{4.14}
\\
\intertext{Furthermore, due to the relations
(\[4.45]) and (\[4.12]) we have}
\ol{(k,l)}\ts\ts\ty_k^{\ts i}\ts\ts\ol{(k,l)}=0\ts;
&\quad
i=1,2,\ts\ldots\ts;\qquad k\neq l.
\Tag{4.145}
\endalign
$$

By definition, the elements $b(\ga)$ where $\ga$ runs through the set
of graphs $\G(n)$, constitute a linear basis in the algebra $B(n,N)$.
Any edge of a graph $\ga\in\G(n)$ of the form
$\{k,l\}$ or $\{\bar k,\bar l\ts\}$ will be called {\it horizontal}. 
If $k<l$ then the vertex $k$ or $\bar k$ will be called the {\it left end}
of the horizontal edge $\{k,l\}$ or $\{\bar k,\bar l\ts\}$ respectively.
The vertex $l$ or $\bar l$ will be then called the {\it right end}.

The number of horizontal edges in a graph $\ga\in\G(n)$ is even.
If this number~is~ $2r$, the element $b(\ga)\in B(n,N)$
has the form 
$
\ts\ol{(k_1,l_1)}\ts\ldots\ts\ol{(k_r,l_r)}\cdot s\ts
$
where $s\in S(n)$ and all $k_1,l_1,\dots,k_r,l_r$ are pairwise distinct.
The elements $b(\ga)$ where the graph $\ga$ has $2r$ horizontal edges
or more, span a two-sided ideal in $B(n,N)$.

\proclaim{Lemma 4.4}
Let $u$ be a monomial in $u_1,\dots,u_n\ts$.
For any two graphs $\ga,\ga\pr\in\G(n)$ we have the equality
in the algebra $\gr\Wen$
$$
b(\ga)\ts u\ts b(\ga\pr)=
\varepsilon\cdot u\pr\ts b(\ga)\ts b(\ga\pr)\ts u\prpr
\Tag{4.15}
$$
where $\varepsilon\in\{\ts1,0,-1\ts\}$
and $u\pr,u\prpr$ are certain monomials in $u_1,\dots,u_n\ts$.
\endproclaim

\demo{Proof}
Let $2r$ and $2r\pr$ be the numbers of horizontal edges in the
graphs $\ga$ and $\ga\pr$ respectively. We will employ the induction
on the minimum of $r,r\pr$ and on the degree of $u$. If each of these
two numbers is zero we have nothing to prove.

Suppose that $r,r\pr\geqslant1$ and $u\neq1$.
By the relations (\[4.12]) we can assume~that
$$
\align
b(\ga)
&=\ol{(k_1,l_1)}\ts\ldots\ts\ol{(k_r,l_r)}=b,
\\
b(\ga\pr)
&=\ol{(k_1^\prime,l^\prime_1)}\ts\ldots\ts\ol{(k^\prime_r,l^\prime_r)}=b\pr
\endalign
$$
where $k_1,l_1,\dots,k_r,l_r$ are pairwise distinct and so are
$k_1^\prime,l^\prime_1,\dots,k^\prime_r,l^\prime_r$.
Consider the monomial
$
u=u_1^{i_1}\ldots u_n^{i_n}\ts.
$
Choose any index $k\in\{1,\dots,n\}$ such that $i_k\neq0$.
If
$$
\align
k\notin\{k_1,l_1,\dots,k_r,l_r\}
&\quad\text{or}\quad
k\notin\{k_1^\prime,l^\prime_1,\dots,k^\prime_r,l^\prime_r\}
\\
\intertext{then respectively} 
b\ts u_k^{i_k}=u_k^{i_k}\ts b
&\quad\text{or}\quad
u_k^{i_k}\ts b\pr=b\pr\ts u_k^{i_k}
\endalign
$$
by (\[4.13]). Then we obtain the equality (\[4.15])
by the inductive assumption.

\smallskip

Now suppose that 
$k=k_j=k^\prime_{j^\prime}$ for some $j$ and $j^\prime$.
Denote $l_j=l$ and $l_{j^\prime}=l^\prime$.
Let~$b\prpr$ denote the product obtained from $b\pr$ by removing the
factor $\ol{(k,l^\prime)}$.
If $\l=l^\prime$ then by the relations (\[4.13]) to (\[4.145]) we have
$$
b\ts u\ts b\pr=(-1)^{i_l}\ts b\ts u_k^{i_k+i_l}\cdot
\prod_{m\neq k,l}u_m^{i_m}\ts\cdot\ol{(k,l)}\ts b\prpr
=
(-1)^{i_l}\ts b\ts u_k^{i_k+i_l}\ts\ol{(k,l)}\cdot
\prod_{m\neq k,l}u_m^{i_m}\ts\cdot b\prpr=0.
$$
Suppose that $\l\neq l^\prime$. Then by the relations (\[4.13]) and (\[4.14])
we have
$$
b\ts u\ts b\pr=(-1)^{i_k}\ts b\ts u_l^{i_k+i_{l^\prime}}\cdot
\prod_{m\neq k,l^\prime}u_m^{i_m}\ts\cdot\ol{(k,l^\prime)}\ts b\prpr
=
(-1)^{i_k}\ts b\ts\ol{(k,l^\prime)}\ts u_l^{i_k+i_{l^\prime}}\cdot
\prod_{m\neq k,l^\prime}u_m^{i_m}\ts\cdot b\prpr.
$$
We have 
$\ts b(\ga\prpr)=b\ts\ol{(k,l^\prime)}\ts$
for a certain graph $\ga\prpr\in\G(n)$. The number of horizontal
edges in the latter graph equals $2r$ as well as in the graph $\ga$.
Thus we obtain the equality (\[4.15]) again
by the inductive assumption\ $\square$
\enddemo

\np
Consider any graph $\ga\in\G(n)$ with exactly $2r$ horizontal edges.
Let
$$
k_1,\dots,k_r\ts,\ts\bar k_1^\prime,\dots\bar k^\prime_r
\quad\text{and}\quad 
l_1,\dots,l_r\ts,\ts\bar l_1^\prime,\dots\bar l^\prime_r
$$
be all the left ends and the right ends of the horizontal edges respectively.

\proclaim{Lemma 4.5}
For any two monomials $u$ and $u\pr$ in $u_1,\dots,u_n$ we have the
equality in the algebra $\gr\Wen$
$$
u\ts b(\ga)\ts u\pr=\varepsilon\cdot
u_1^{i_1}\ldots\ts u_n^{i_n}
\ts\ts b(\ga)\ts\ts
u_1^{j_1}\ldots\ts u_n^{j_n}
$$
where $\varepsilon\in\{\ts1,0,-1\ts\}$ and
$$
k\in\{l_1,\dots,l_r\}\ts\ts\Rightarrow\ts\ts i_k=0\ts;
\ \quad
j_k\neq0\ts\ts\Rightarrow\ts\ts k\in\{l^\prime_1,\dots,l^\prime_r\}\ts.
\Tag{4.16}
$$
\endproclaim

\demo{Proof}
Due to the relations (\[4.12]) we can assume that
$b(\ga)=\ol{(k_1,l_1)}\ts\ldots\ts\ol{(k_r,l_r)}\ts,$
$$
k_1=k_1^\prime,\dots,k_r=k_r^\prime
\quad\text{and}\quad
l_1=l_1^\prime,\dots,l_r=l_r^\prime\ts.
$$
The required statement follows then directly from the relations
(\[4.13]) and (\[4.14])\ $\square$
\enddemo

\np
Any product in the algebra $\Wen$ of the form
$$
y_1^{i_1}\ldots\ts y_n^{i_n}
\ts\ts b(\ga)\ts\ts
y_1^{j_1}\ldots\ts y_n^{j_n}
\ts\cdot\ts
w_2^{h_2}\ts w_4^{h_4}\dots
\Tag{4.17}
$$
will be called a {\it regular monomial}
if the exponents $i_1,\dots,i_n$ and $j_1,\dots,j_n$ satisfy
the conditions
(\[4.16]). The two theorems below are the main results
of this section.

\proclaim{Theorem 4.6}
All the regular monomials {\rm (\[4.17])}
constitute a basis in $\Wen$.
\endproclaim

\np
By the relations (\[2.2]) to (\[2.4]) and (\[2.45]) for every
$m=0,1,2,\ts\ldots$ the assignments
$$
s_k\mapsto s_{m+k},
\quad
\bs_k\mapsto \bs_{m+k},
\quad
y_k\mapsto x_{m+k},
\quad
w_i\mapsto z_{m+1}^{(i)}
$$
define a homomorphism
$$
\pi_m:\Wen\to B(m+n,N).
$$
The homomorphism $\pi_0$ coincides with (\[4.455555]).
Furthermore, by Lemma~2.1 the image of the homomorphism
$\pi_m$ commutes with the subalgebra $B(m,N)$ in $B(m+n,N)$. 

\proclaim{Theorem 4.7}
The kernels of $\pi_0,\pi_1,\pi_2,\ts\ldots$
have zero intersection.
\endproclaim 

\np
Due to Lemmas 4.4\ts,\ts4.5 and to the equalities (\[4.475])
every element of the algebra $\Wen$ can be expressed as a linear
combination of regular monomials. Thus both Theorems 4.6 and 4.7
follow from the next lemma; cf. [O,\ts Lemma~2.1.11].
Fix any finite set $\Cal F$ of regular monomials (\[4.17]).
Let $m$ be the maximum of the sums
$$
i_1+j_1+\ldots+i_n+j_n+2\ts h_2+4\ts h_4+\ldots
$$
corresponding to the monomials from the set $\Cal F$. 
Consider any linear combination $F\in\Wen$ of some monomials
from $\Cal F$ with non-zero coefficients.

\proclaim{Lemma 4.8}
Suppose that $\pi_m(F)=0$. Then monomials on which
\text{the maximum} $m$ is attained, do not appear in $F$.
\endproclaim

\demo{Proof}
Due to Corollary~4.3 
we have the equality 
$$
Q_{m+1}(u)=
\operatorname{exp}
\sum_{i=1,3,\dots}2\ts
\Bigl(
(N-1)^i/2^i+\sum_{j=1}^m\ts(x_j+1)^i-2x_j^{\ts i}+(x_j-1)^i
\ts\Bigr)
\ts{u^{-i}}/{i}
$$
of formal power series in $u^{-1}$. Therefore for each $i=2,4,\ts\ldots$
by the equality
(\[2.4666666])
the element $z_{m+1}^{(i)}\in B(m,N)$ is a symmetric polynomial in
$x_1,\ldots,x_m$ of the form
$$
2i\ts\bigl(x_1^{i-1}+\ldots+x_m^{i-1}\bigr)
+\text{terms of smaller degrees.}
$$

Consider the subset $\Cal F^{\,\prime}\subset\Cal F$ formed by all
the monomials where the maximum $m$ is attained. Then consider
the subset $\Cal F^{\,\prime\prime}\subset\Cal F^{\,\prime}$ formed
by the monomials with the minimal number of horizontal edges 
in the corresponding graphs $\ga$. Let $2r$ be that minimum.
It suffices to prove that the monomials from $\Cal F^{\,\prime\prime}$
do not appear~in~$F$.

Choose any regular monomial (\[4.17]) from the set $\Cal F$.
The image of this monomial with respect to the homomorphism $\pi_m$ is
a certain linear combination $f$ of the elements $b(\Ga)\in B(m+n,N)$
where $\Ga\in\G(m+n)$.
Denote by $\G$ be the subset in $\G(m+n)$ consisting of those graphs
which have:

\smallskip

\np
\quad- exactly $2r$ horizontal edges;

\np
\quad- no vertical edges of the form $\{k,\bar k\ts\}$
where $k\leqslant m\ts$;

\np
\quad- no horizontal edges of the form $\{k,l\}$ or $\{\bar k,\bar l\ts\}$
where $k,l\leqslant m\ts$.

\smallskip
\np
Consider the terms of $f$ corresponding
to the graphs $\Ga\in\G$. Such terms appear
in $f$ only if the chosen monomial
belongs to the subset $\Cal F^{\,\prime\prime}$. Suppose this is the case.
Then amongst those terms are the products
$$
\align
&\prod_{k=1}^n\ts\ts
(m_{k1}\ts,\ts m+k)
\ldots
(m_{k\ts i_k}\ts,\ts m+k) 
\ts\cdot\ts\pi_m\bigl(b(\ga)\bigr)\ts\ts\times
\Tag{4.20}
\\
\qquad\ \qquad&\prod_{k=1}^n\ts\ts
(m_{k1}^{\ts\prime}\ts,\ts m+k)
\ldots
(m_{k\ts j_k}^{\ts\prime}\ts,\ts m+k) 
\cdot
\prod_{i=2,4,\ldots}\ts\ts\ts
\prod_{j=1}^{h_i}\ts\ts
2\ts(m_{ij1}\ts,\ldots,\ts m_{iji})
\endalign
$$
where the juxtaposition of the sequences
$$
m_{k1}
\ts,\ts\ldots,\ts
m_{k\ts i_k}
\ts,
m_{k1}^{\ts\prime}
\ts,\ts\ldots,\ts
m_{k\ts j_k}^{\ts\prime}
\ts;
\qquad k=1,\ts\ldots\ts,n
$$
and
$$
m_{ij1}\ts,\ldots,\ts m_{iji}\ts;
\qquad
j=1,\ts\ldots\ts,h_i\ts;
\quad
i=2,4,\ldots
$$
runs through the set of all permutations
of the sequence $1,2,\ldots,m$.
All these terms will be called the {\it leading terms} of $f$. 
Note that the parameters
$$
\ga\ts,\ts i_1,\dots,i_n\ts,\ts j_1,\dots,j_n\ts,\ts h_2,h_4,\dots
$$
can be uniquely restored from any of these leading terms.

All the non-leading terms of $f$ corresponding
to graphs $\Ga\in\G$ can be obtained from 
the products (\[4.20])
by certain non-empty 
sets of the following replacements. One can replace 
the factor in (\[4.20])
$$
(m_{k1}\ts,\ts m+k)
\ldots
(m_{k\ts i_k}\ts,\ts m+k) 
\quad\text{by}\quad
\ol{(m_{k1}\ts,\ts m+k)}
\ldots
\ol{(m_{k\ts i_k}\ts,\ts m+k)}\cdot(-1)^{i_k} 
\nopagebreak
$$
provided the vertex $k$ of the graph $\ga$ is the left end of a
horizontal edge. One can also replace any factor 
$$
(m_{k1}^{\ts\prime}\ts,\ts m+k)
\ldots
(m_{k\ts j_k}^{\ts\prime}\ts,\ts m+k) 
\quad\text{by}\quad
\ol{(m_{k1}^{\ts\prime}\ts,\ts m+k)}
\ldots
\ol{(m_{k\ts j_k}^{\ts\prime}\ts,\ts m+k)}\cdot(-1)^{j_k}. 
$$
Due to the conditions (\[4.16])
the terms so obtained are not proportional to any leading term
in the image with respect to $\pi_m$ of any monomial from
$\Cal F^{\,\prime\prime}$. 
This observation completes the proof
\ $\square$
\enddemo

\np
We will now compare the algebra $\Wen$ with the
{\it degenerate affine Hecke algebra} $\Hen$
from [C1,C2] and [D].
The latter algebra is generated by the group algebra $\CC[S(n)]$
and the pairwise commuting elements $v_1,\dots,v_n$ 
subjected to the relations
$$
\align
s_k\ts v_l=v_l\ts s_k\ts,
&\quad
l\neq k,k+1;
\\
s_k\ts v_k-v_{k+1}\ts s_k=-1\ts,
&\quad
s_k\ts v_{k+1}-v_k\ts s_k=1\ts.
\endalign
$$
By the relations (\[4.2]) to (\[4.45]) we have the following
corollary to Theorem~4.6.

\proclaim{Corollary 4.9}
For any $f_2,f_4,\ldots\in\CC$ the maps
$s_k\mapsto s_k$,
$\bs_k\mapsto0$,
$y_k\mapsto v_{k}$ and
$$
w_i\mapsto f_i,
\quad
i=2,4,\dots
$$
determine a homo\-morphism of the algebra $\Wen$ onto $\Hen$.
\endproclaim

\np
The subalgebra in $\Hen$ generated by the elements
$v_1,\dots,v_n$ is maximal commutative. The centre of the algebra $\Hen$
consists of all symmetric polynomials in $v_1,\dots,v_n\ts$.
For the proofs of these two statements see [C2\ts,\ts Section 1].
The next corollary provides analogues of these statements
for the algebra $\Wen$.

\proclaim{Corollary 4.10} 
The subalgebra in $\Wen$ generated by the elements
$y_1,\dots,y_n$ and $w_1,w_2,\ts\ldots$ is maximal commutative.
The elements $y_1^i+\ldots+y_n^i$ with $i=1,3,\ldots$
and $w_i$ with $i=2,4,\ldots$
generate the centre of the algebra $\Wen$.
\endproclaim 

\bigskip

\heading
References
\endheading

\itemitem{[Br]}
{R.\, Brauer,}
On algebras which are connected with the semisimple continuous groups,
{\it Ann.\, of Math.}
{\bf 38}
(1937),
857--872.

\itemitem{[BW]}
{J.\, S.\, Birman and H.\, Wenzl,}
Braids, link polynomials and a new algebra,
{\it Trans. Amer.\, Math.\, Soc.}
{\bf 313}
(1989),
249--273.

\itemitem{[B1]}
{W.\, P.\, Brown,}
An algebra related to the orthogonal group,
{\it Michigan Math.\, J.}
{\bf 3}
(1955--1956),
1--22.

\itemitem{[B2]}
{W.\, P.\, Brown,}
The semisimplicity of $\omega_f^n\ts$,
{\it Ann.\, of Math.}
{\bf 63}
(1956),
324--335.

\itemitem{[C1]}
{I.\,V.\,Cherednik},
{On special bases of irreducible finite-dimensional representations of the
degenerate affine Hecke algebra},
{\it Funct. Analysis Appl.}
{\bf 20}
(1986),
76--78.

\itemitem{[C2]}
{I.\,V.\,Cherednik},
{A unification of Knizhnik-Zamolodchikov and Dunkl operators
via affine Hecke algebras},
{\it Invent. Math.}
{\bf 106}
(1991),
411--431.

\itemitem{[D]}
{V.\,G.\,Drinfeld},
{Degenerate affine Hecke algebras and Yangians},
{\it Funct. Analysis Appl.}
{\bf 20}
(1986),
56--58.

\itemitem{[EK]}
{N.\, El Samra and R.\, C.\, King,}
Dimensions of irreducible representations of the classical Lie groups,
{\it J.\, Phys.}
{\bf A12}
(1979),
2317--2328.

\itemitem{[HW]}
{P.\, Hanlon and D.\, Wales,}
A tower construction for the radical in Brauer's
centralizer algebras,
{\it J.\, Algebra}
{\bf 164}
(1994),
773--830.

\itemitem{[JMO]}
{M.\, Jimbo, T.\, Miwa and M.\, Okado,}
Solvable lattice models related to the vector representations of the
classical Lie algebras,
{\it Comm.\, Math.\, Phys.} 
{\bf 116}
(1988),
507--525.

\itemitem{[Ju]}
{A.\,-A.\, A.\, Jucys,}
Symmetric polynomials and the center of the symmetric group ring,
{\it Rep.\, Math.\, Phys.}
{\bf 5}
(1974),
107--112.

\itemitem{[Ke]}
{S.V. Kerov,}
Realization of representations of the Brauer semigroup,
{\it J.\,Soviet Math.}
{\bf 47}
(1989),
2503--2507.

\itemitem{[Ki]}
{R.\, C.\, King,}
Modification rules and products of irreducible representations of the
unitary, orthogonal, and symplectic groups,
{\it J.\, Math.\, Phys.}
{\bf 12}
(1971),
1588--1598.

\itemitem{[L]}
{D.\, E.\, Littlewood,}
Products and plethysms of characters with orthogonal, symplectic and
symmetric groups,
{\it Canad.\, J.\, Math.}
{\bf 10}
(1958),
17--32.

\itemitem{[Ma]}
{I.\, G.\, Macdonald,}
\lq\lq\ts Symmetric Functions and Hall Polynomials\ts\rq\rq,
Clarendon Press,
Oxford,
1979.

\itemitem{[Mk]}
{J.\, Murakami,}
The representations of the $q$-analogue of Brauer's centralizer
algebras and the Kauffman polynomial of links,
{\it Publ.\, RIMS} 
{\bf 26}
(1990),
{935--945}.

\itemitem{[Mo]}
{V.\, F.\, Molchanov,}
On matrix elements of irreducible representations of symmetric group,
{\it Vestn.\, Mosk.\, Univ.}
{\bf 1}
(1966),
52--57.

\itemitem{[Mp]}
{G.\, E.\, Murphy,}
A new construction of Young's seminormal representation of the symmetric
group,
{\it J.\, Algebra}
{\bf 69}
(1981),
287--291.

\itemitem{[O]}
{G.\,I.\,Olshanski\u\i},
Representations of infinite-dimensional classical groups,
limits of enveloping algebras, and Yangians,
{\it in} \lq\lq\ts Topics in Representation~Theory\ts\rq\rq,
{\it Adv.\,Soviet Math.}
{\bf 2},
Amer.\,Math.\,Soc.,
Providence,
1991,
pp. 1--66.

\itemitem{[RW]}
{A.\, Ram and H.\, Wenzl,}
Matrix units for centralizer algebras,
{\it J.\, Algebra}
{\bf 145}
(1992),
378--395.

\itemitem{[V1]}
{A.\, M.\, Vershik,}
Local algebras and a new version of Young's orthogonal form,
{\it in} \lq\lq\ts Topics in Algebra\ts\rq\rq,
{\it Banach Center Publications}
{\bf 26},
Polish Scientific,
Warsaw,
1990,
pp. 467--473.

\itemitem{[V2]}
{A.\, M.\, Vershik,}
Local stationary algebras,
{\it Amer.\, Math.\, Soc.\, Translations} 
{\bf 148}
(1991),
1--13.

\itemitem{[W]}
{H.\, Wenzl,}
On the structure of Brauer's centralizer algebras,
{\it Ann.\, of Math.}
{\bf 128}
(1988),
173--193.

\itemitem{[Wy]}
{H.\, Weyl,}
\lq\lq\ts Classical groups, their invariants and representations\ts\rq\rq,
Pri\-nce\-ton University Press,
Princeton NJ,
1946.

\itemitem{[Y]}
{A.\, Young,}
On quantitative substitutional analysis VI,
{\it Proc.\,London Math. Soc.}
{\bf 31}
(1931),
253--289.

\bye